\documentstyle[12pt, amsfonts]{amsart}
\begin{document}
\def\R{{\mathbb R}}
\def\Z{{\mathbb Z}}
\def\C{{\mathbb C}}
\newcommand{\trace}{\rm trace}
\newcommand{\Ex}{{\mathbb{E}}}
\newcommand{\Prob}{{\mathbb{P}}}
\newcommand{\E}{{\cal E}}
\newcommand{\F}{{\cal F}}
\newtheorem{df}{Definition}
\newtheorem{theorem}{Theorem}
\newtheorem{lemma}{Lemma}
\newtheorem{pr}{Proposition}
\newtheorem{co}{Corollary}
\def\n{\nu}
\def\sign{\mbox{ sign }}
\def\a{\alpha}
\def\N{{\mathbb N}}
\def\A{{\cal A}}
\def\L{{\cal L}}
\def\X{{\cal X}}
\def\F{{\cal F}}
\def\c{\bar{c}}
\def\v{\nu}
\def\d{\delta}
\def\diam{\mbox{\rm dim}}
\def\vol{\mbox{\rm Vol}}
\def\b{\beta}
\def\t{\theta}
\def\l{\lambda}
\def\e{\varepsilon}
\def\colon{{:}\;}
\def\pf{\noindent {\bf Proof :  \  }}
\def\endpf{ \begin{flushright}
$ \Box $ \\
\end{flushright}}

\title[Positive definite functions]{Positive definite functions and  
multidimensional versions of random variables}

\author{Alexander Koldobsky}

\address{Department of Mathematics\\ 
University of Missouri\\
Columbia, MO 65211}

\email{koldobsk@@math.missouri.edu}

\begin{abstract} 
We say that a random vector $X=(X_1,...,X_n)$ in $\R^n$ is an $n$-dimensional version
of a random variable $Y$ if for any $a\in \R^n$ the random variables $\sum a_iX_i$
and $\gamma(a) Y$ are identically distributed, where $\gamma:\R^n\to [0,\infty)$ is 
called the standard of $X.$ An old 
problem is to characterize those functions $\gamma$ that can appear as the standard
of an $n$-dimensional version. In this paper, we prove the conjecture of Lisitsky that
every standard must be the norm of  a space that embeds in $L_0.$ This result 
is almost optimal, as the norm of any finite dimensional subspace of $L_p$ with $p\in (0,2]$
is the standard of an $n$-dimensional version ($p$-stable random vector) by the 
classical result of P.L\`evy. An equivalent formulation is that if  a function of the form 
$f(\|\cdot\|_K)$ is positive definite on $\R^n,$ where $K$ is an origin symmetric star 
body in $\R^n$ and $f:\R\to \R$ is an even continuous function, then either the space  
$(\R^n,\|\cdot\|_K)$ embeds in $L_0$ or $f$ is a constant function.  Combined with known 
facts about embedding in $L_0,$ this result leads to several generalizations of the solution 
of Schoenberg's problem on positive definite functions.
\end{abstract}  
\maketitle

\section{Introduction}

Following Eaton [E], we say that a random vector $X=(X_1,...,X_n)$ 
is an {\it $n$-dimensional version} 
of a random variable $Y$ if there exists a function $\gamma: \R^n\to \R,$ called 
the {\it standard} of $X,$ such that  $\gamma(a)>0$ for every $a\in \R^n,\ a\neq 0,$ and for 
every $a\in \R^n$ the random variables 
\begin{equation} \label{vers}
\sum_{i=1}^n a_i X_i \quad {\rm and} \quad \gamma(a) Y
\end{equation}
are identically distributed. We assume that $n\ge 2$ and $P\{Y=0\}<1.$  A problem posed 
by Eaton is to characterize all $n$-dimensional versions,
and, in particular, characterize all functions $\gamma$ that can appear 
as the standard of  an $n$-dimensional version. 

It is easily seen [M3], [Ku] that every standard $\gamma$ is an even homogeneous of degree 
1 non-negative (and equal to zero only 
at zero) continuous function on $\R^n.$  This means that $\gamma=\|\cdot\|_K$
is the Minkowski functional of some origin symmetric star body $K$ in $\R^n.$
Recall that a closed bounded set $K$ in $\R^n$ is called {\it a star body} if 
every straight line passing through the origin crosses the boundary of $K$ 
at exactly two points, the origin is an interior point of $K$
and {\it the Minkowski functional} of $K$ defined by
$\|x\|_K =\min\{s\ge 0:\ x\in sK\}$
is a continuous function on $\R^n.$ Note that the class of star bodies includes 
convex bodies containing the origin in their interior.

Eaton [E] proved that a random vector is an $n$-dimensional version with the standard $\|\cdot\|_K$
if and only if its characteristic functional has the form $f(\|\cdot\|_K),$ where $K$ is an
origin symmetric star body in $\R^n$ and $f$ is an even continuous non-constant function on $\R$
(see also [K3, Lemma 6.1] ).  By Bochner's theorem,
this means that the function $f(\|\cdot\|_K)$ is positive definite.
Recall that a complex valued function $f$ defined on $\R^n$ is called {\it positive definite} 
on $\R^n$ if, for every finite sequence $\{x_i\}_{i=1}^m$ in $\R^n$ and every choice
of complex numbers $\{c_i\}_{i=1}^m$, we have
$$\sum_{i=1}^m\sum_{j=1}^mc_i\bar{c}_jf(x_i-x_j)\ge 0.$$

Thus, Eaton's problem is equivalent to characterizing the classes $\Phi(K)$
consisting of even continuous functions $f:\R\to \R$ for which $f(\|\cdot\|_K)$ is
a positive definite function on $\R^n.$ In particular, $\|\cdot\|_K$ appears as the standard
of an $n$-dimensional version if and only if the class $\Phi(K)$ is non-trivial, i.e.
contains at least one non-constant function. In some places throughout the paper we write $\Phi(E_K)$
instead of $\Phi(K),$ where $E_K=(\R^n,\|\cdot\|_K)$ is the space whose unit ball is $K.$

The problem of characterization of positive definite norm dependent functions has a long
history and goes back to the work of L\`evy and Schoenberg in the 1930s.
L\`evy [Le] proved that, for any finite dimensional subspace $(\R^n,\|\cdot\|)$ of $L_q$
with $0<q\le 2,$ the function $g=\exp(-\|\cdot\|^q)$ is positive definite on $\R^n$, 
and any random vector $X=(X_1,...,X_n)$ in $\R^n$, whose characteristic functional is $g$, 
satisfies the property  (\ref{vers}). This result gave a start to the theory of stable processes that
has numerous applications to different areas of mathematics. The concept of an $n$-dimensional
version is a generalization of stable random vectors.

In 1938, Schoenberg [S1,S2] found a connection between positive definite functions 
and the embedding theory of metric spaces. In particular, Schoenberg [S1] posed the 
problem of finding the exponents $0<p\le 2$ for which the function $\exp(-\|\cdot\|_q^p)$ is positive 
definite on $\R^n$, where 
$$\|x\|_q= \left(|x_1|^q+...+|x_n|^q\right)^{1/q}$$
is the norm the space $\ell_q^n$ with $2<q\leq \infty.$ 
This problem had been open for more than fifty years.
For $q=\infty,$ the problem was solved in 1989 by 
Misiewicz [M2], and for $2<q<\infty,$ the answer was 
given in [K1] in 1991 (note that, for $1\le p\le 2,$ Schoenberg's question was answered 
earlier by Dor [D], and the case $n=2,\ 0<p\le 1$ was established in [F], [H], [L]). 
The answers turned out to be the same in both cases: the function
$\exp(-\|\cdot\|_q^p)$ is not positive definite for any $p\in (0,2]$ if $n\ge 3$, and for $n=2$ the function 
is positive definite if and only if $0<p\leq 1.$  Different and independent 
proofs of Schoenberg's problems were given by Lisitsky [Li1] and  Zastavnyi [Z1, Z2] 
shortly after the paper [K1] appeared. For generalizations of the solution of Schoenberg's
problem, see [KL].

The solution of Schoenberg's problem can be interpreted in terms of isometric 
embeddings of normed spaces. In fact, the result of Bretagnolle, Dacunha-Castelle
and Krivine [BDK] shows that a normed space embeds isometrically in $L_p,\ 0<p\le 2$
if and only if the function $\exp(-\|\cdot\|^p)$ is positive definite. Hence, the answer
to Schoenberg's problem means that that the spaces $\ell_q^n,\ q>2,\ n\ge 3$ do not
embed isometrically in $L_p$ with $0<p\le 2.$ 

The classes $\Phi(K)$ have been studied by a number
of authors. Schoenberg [S2] proved that
$f\in \Phi(\ell_2^n)$ if and only if 
$$ f(t)=\int_0^{\infty} \Omega_n(tr)\ d\lambda(r)$$
where $\Omega_n(|\cdot|_2)$ is the Fourier transform of 
the uniform probability measure on the sphere $S^{n-1},$
$|\cdot|_2$ is the Euclidean norm in $\R^n,$
and $\lambda$ is a finite measure on $[0,\infty).$ 
In the same paper, Schoenberg proved an infinite dimensional 
version of this result: $f\in \Phi(\ell_2)$ if and only if
$$ f(t)=\int_0^{\infty} \exp(-t^2r^2)\ d\lambda(r).$$
Bretagnolle, Dacunha-Castelle and Krivine [BDK] proved
a similar result for the classes $\Phi(\ell_q)$ for 
all $q\in (0,2)$ (one just has to replace 2 by $q$ in the formula),
and showed that for $q>2$ the classes $\Phi(\ell_q)$ (corresponding
to infinite dimensional $\ell_q$-spaces) are trivial,
i.e. contain constant functions only. Cambanis, Keener and Simons 
[CKS] obtained a similar representation for the classes
$\Phi(\ell_1^n).$ Richards [R] and Gneiting [G] partially characterized the 
classes $\Phi(\ell_q^n)$ for $0<q<2.$ Aharoni, Maurey and Mityagin
[AMM] proved that if $E$ is an infinite dimensional Banach space 
with a symmetric basis $\{e_n\}_{n=1}^{\infty}$ such that 
$$\lim \limits_{n\to \infty} \frac{\|e_1+\dots+e_n\|}{n^{1/2}}=0,$$
then the class $\Phi(E)$ is trivial.  Misiewicz [M2] proved that for $n\geq 3$ the
classes $\Phi(\ell_{\infty}^n)$ are trivial, and Lisitsky [Li1] and Zastavnyi [Z1], [Z2] 
showed the same for the classes $\Phi(\ell_q^n),\ q>2,\ n\ge 3.$ One can find more related 
results and references  in [M3], [K3].

In all the results mentioned above the
classes $\Phi(K)$ appear to be non-trivial only if
$K$ is the unit ball of a subspace of $L_q$ with $0<q\leq 2.$
An old conjecture, explicitly formulated for the first time by
Misiewicz [M1],  is that the class $\Phi(K)$ can be non-trivial only in this case.
A slightly weaker conjecture was formulated by Lisitsky [Li2]:  if  the class 
$\Phi(K)$ is non-trivial,  then the space $(\R^n,\|\cdot\|_K)$
embeds in $L_0.$  The concept of embedding in $L_0$ was introduced
and studied in [KKYY], the original conjecture of Lisitsky was 
in terms of the representation (\ref{logrepr}):
\begin{df} We say that a space $(\R^n, \|\cdot\|_K)$ 
embeds in $L_0$ if there exist a finite Borel
measure $\mu$ on the sphere $S^{n-1}$ and a constant $C\in \R$ so that, for every
$x\in\R^n$,
\begin{equation} \label{logrepr}
\ln \|x\|_K =\int_{S^{n-1}} \ln |(x, \xi ) |\  d\mu(\xi) + C.
\end{equation}
\end{df}
It is quite easy to confirm the conjectures of Misiewicz and Lisitsky under 
additional assumptions that $f$ or its Fourier transform have finite moments of certain
orders; see [Mi1], [Ku], [Li2], [K4]. 

\medbreak
In this article we prove the conjecture of Lisitsky in its full strength:

\begin{theorem} \label{main} Let $K$ be an origin symmetric star body in $\R^n,\ n\ge 2$ and suppose that
there exists an even non-constant continuous function $f:\R \mapsto \R$ such that  $f(\|\cdot\|_K)$
is a positive definite function on $\R^n.$ Then the space $(\R^n,\|\cdot\|_K)$ embeds 
in $L_0.$
\end{theorem}

\begin{co} If a function $\gamma$ is the standard of an $n$-dimensional version of a random
variable, then there exists an origin symmetric star body $K$ in $\R^n$ such that $\gamma=\|\cdot\|_K$
and the space $(\R^n,\|\cdot\|_K)$ embeds in $L_0.$
\end{co}

In the last section of the paper we use known results about embedding in $L_0$ to
point out rather general classes of normed spaces for which the classes $\Phi$ are trivial
and whose norms cannot serve as the standard of an $n$-dimensional version.

\section{Proof of Theorem \ref{main} }

As usual, we denote by ${\mathcal S}(\R^n)$ the space of infinitely differentiable rapidly
decreasing functions on $\R^n$ (Schwartz test functions), and by ${\mathcal S}^{'}(\R^n)$ the space
of distributions over ${\mathcal S}(\R^n).$ If $\phi\in {\mathcal S}(\R^n)$ and $f\in {\mathcal S}^{'}(\R^n)$ 
is a locally integrable function with power growth at infinity, then the action of $f$ on $\phi$ is defined by
$$\langle f, \phi \rangle = \int_{\R^n} f(x)\phi(x)\ dx.$$

We say that a distribution is positive (negative) outside of the origin in $\R^n$ if it assumes
non-negative (non-positive) values on non-negative test functions with compact support
outside of the origin.

The Fourier transform of a distribution $f$ is defined by
$\langle\hat{f}, {\phi}\rangle= \langle f, \hat{\phi} \rangle$
for every test function $\phi.$ A distribution is positive definite
if its Fourier transform is a positive distribution.
\medbreak
We use the following Fourier analytic characterization of embedding in $L_0$ proved 
in [KKYY, Th.3.1]:
\begin{pr} \label{log} Let $K$ be an origin symmetric star body in $\R^n$.
The space $(\R^n, \|\cdot\|_K)$ embeds in $L_0$ if and only if the Fourier transform of
$\ln \|x\|_K$ is a negative distribution outside of the origin in $\R^n$.
\end{pr}
\medbreak

Now we are ready to start the proof of Theorem \ref{main}.
\bigbreak
{\bf Proof of Theorem \ref{main}.} We write $\|\cdot\|$ instead of $\|\cdot\|_K.$ By Bochner's theorem, the function
$f(\|\cdot\|)$ is the Fourier transform of a finite measure $\mu$ on $\R^n.$ We can 
assume that $f(0)=1,$ and, correspondingly, $\mu$ is a probability measure. 
The function $f$ is positive definite on $\R,$ as the
restriction of a positive definite function, therefore, $|f(t)|\le f(0)=1$ for every
$t\in \R$ (see [VTC, p.188]).

Let $\phi$ be an even non-negative test function supported outside 
of the origin in $\R^n.$
For every fixed $t>0,$ the function $f(t\|\cdot\|)$ is positive definite on $\R^n,$ so
\begin{equation} \label{fposdef}
\int_{\R} f(t\|x\|)  \hat\phi(x)\ dx  = \big\langle \left(f(t\|\cdot\|)\right)^\wedge, \phi(x) \big\rangle \ge 0.
\end{equation}

For any $\e\in (0,1),$ the integral
\begin{equation} \label{expr}
g(\e)= \int_{\R^n} \left( \int_0^1 t^{-1+\e} f(t\|x\|)dt + \int_1^\infty t^{-1-\e} f(t\|x\|) dt \right) \hat\phi(x)dx 
\end{equation}
converges absolutely, because $f$ is bounded by 1 and the function in parentheses is bounded by $2/ \e.$
By the Fubini theorem,  
$$g(\e) = \int_0^1 t^{-1+\e} \left(\int_{\R^n} f(t\|x\|)\hat\phi(x) dx \right) dt  $$
$$+ \int_1^\infty t^{-1-\e} \left( \int_{\R^n} f(t\|x\|) \hat\phi(x)dx \right) dt,$$
so by (\ref{fposdef}) the function $g$ is non-negative:
\begin{equation} \label{nonneg}
g(\e)\ge 0 \quad  {\rm for\ every} \quad \e\in (0,1).
\end{equation}

Now we study the behavior of the function $g$, as $\e\to 0.$ We have
$$g(\e) = \int_{\R^n} \left( \|x\|^{-\e} \int_0^{\|x\|} t^{-1+\e} f(t)\ dt + \|x\|^\e \int_{\|x\|}^\infty t^{-1-\e} f(t)\ dt \right) \hat\phi(x) dx$$
\begin{equation} \label{int1}
= \int_{\R^n} \frac{\|x\|^{-\e}-1}{\e} \e\left( \int_0^{\|x\|} t^{-1+\e} f(t)\ dt \right) \hat\phi(x) dx
\end{equation}
\begin{equation} \label{int2}
 + \int_{\R^n} \frac{\|x\|^\e-1}{\e}\e \left(\int_{\|x\|}^\infty t^{-1-\e} f(t)\ dt \right) \hat\phi(x) dx
 \end{equation}
\begin{equation} \label{int3}
+\int_{\R^n} \left( \int_0^{\|x\|} t^{-1+\e} f(t)\ dt +  \int_{\|x\|}^\infty t^{-1-\e} f(t)\ dt \right) \hat\phi(x) dx.
\end{equation}
We write
$$g(\e) = u(\e) + v(\e) + w(\e),$$
where $u,v,w$ are the functions defined by (\ref{int1}), (\ref{int2}) and (\ref{int3}), respectively.
\bigbreak
We start with the function $w.$

\begin{lemma} \label{l1}
$$\lim_{\e\to 0}\ w(\e) = 0.$$
\end{lemma}

\pf  We can assume that $\e<1/2.$ Fix $a>0.$  Since $\phi$ is supported outside of the origin, we have $\int_{\R^n} \hat\phi(x)dx =0$ and 
$$\int_{\R^n} \left( \int_0^{a} t^{-1+\e} f(t)\ dt +  \int_{a}^\infty t^{-1-\e} f(t)\ dt \right) \hat\phi(x) dx = 0,$$
because the expression in parentheses is a constant.
Subtracting this from (\ref{int3}) we get
$$w(\e) = \int_{\R^n} \left( \int_a^{\|x\|}\left( t^{-1+\e} -t^{-1-\e}\right) f(t)\ dt \right) \hat\phi(x) dx.$$
Now for some $\theta(t,\e)\in [0,2\e],$ 
$$t^{-1-\e}|t^{2\e}-1| = 2\e\ t^{-1-\e}t^{\theta(t,\e)} |\ln t|$$
$$ \le 2\e (1+a^{-3/2}+\|x\|^{-3/2})(|\ln a|+\big|\ln \|x\|\big|),$$
so 
\begin{equation} \label{estw}
|w(\e)|\le 2\e \int_{\R^n} | \|x\|-a|\ (1+a^{-3/2}+\|x\|^{-3/2})(|\ln a|+\big|\ln \|x\|\big|)|\hat\phi(x)|dx.
\end{equation}
By the definition of a star body, $K$ is bounded and contains a Euclidean ball with center at the origin,
so there exist constants $c,d>0$ so that for every $x\in \R^n$
\begin{equation} \label{eucl}
c|x|_2\le \|x\|\le d|x|_2,
\end{equation}
where $|\cdot|_2$ is the Euclidean norm in $\R^n.$
Note that $n\ge 2$ so $|\cdot|_2^{-3/2}$ is a locally integrable function on $\R^n,\ n\ge 2.$ 
Also $\hat\phi$ is a test function and decreases at infinity faster than any power of the 
Euclidean norm. These facts, in conjunction with (\ref{eucl}),  imply that the integral in the right-hand side 
of (\ref{estw}) converges, which proves the lemma. \qed
\bigbreak
We need the following elementary and well known fact. 
\begin{lemma} \label{eps} Let $h$ be a bounded integrable continuous at 0 function on $[0,A],\ A>0.$ Then
$$\lim_{\e\to 0} \e \int_0^A t^{-1+\e} h(t) dt = h(0).$$
\end{lemma}

\pf We can assume that  $\e<1.$ We have 
$$\e \int_0^A t^{-1+\e} h(t) dt $$
$$=  \e \int_0^{{\e}} t^{-1+\e} (h(t)-h(0)) dt +
\e h(0) \int_0^{{\e}} t^{-1+\e} dt + \e \int_{{\e}}^A t^{-1+\e} h(t) dt.$$
The first summand is less or equal to  
$${\e}^{\e} \max_{ t \in [0, \e]} |h(t)-h(0)| \to 0, \quad {\rm as} \quad \e \to 0,$$
because $h$ is continuous at 0.
The second summand is equal to 
$$h(0) \e^{\e}\to h(0),\quad {\rm as} \quad \e \to 0.$$ 
The third summand is less or equal to 
$$|A^\e-\e^\e| \max_{t\in [0,A]} |h(t)| \to 0, \quad \rm{as}\quad \e \to 0. \qed$$

\bigbreak
Now we compute the limit at infinity of the function
$$u(\e)= \int_{\R^n} \frac{\|x\|^{-\e}-1}{\e} \e\left( \int_0^{\|x\|} t^{-1+\e} f(t)\ dt \right) \hat\phi(x) dx.$$

\begin{lemma} \label{l2}
$$\lim_{\e\to 0}\ u(\e) = -f(0) \int_{\R^n}\ln \|x\| \hat\phi(x) dx.$$
\end{lemma}

\pf Using the estimates
$$\Big|\frac{\|x\|^{-\e}-1}{\e}\Big| =\Big|\frac{1}{\e}\int_0^{\e} \|x\|^{-\theta}\ln \|x\| d\theta\Big|\le \big|\ln\|x\|\big|(1+\|x\|^{-1})$$
and 
$$\Big|\e \int_0^{\|x\|} t^{-1+\e} f(t) dt \Big| \le  \|x\|^{\e}\le \|x\|+1,$$
we see that the functions under the integral over $\R^n$ in $u(\e)$
are dominated by an integrable function 
$$\big|\ln\|x\|\big|(1+\|x\|^{-1})(\|x\|+1) |\hat\phi(x)|$$
of the variable $x$ on $\R^n.$ Clearly, for $x\neq 0,$
$$\lim_{\e\to 0} \frac{\|x\|^{-\e}-1}{\e} = -\ln \|x\|.$$
Also, by Lemma \ref{eps}, for every $x\in \R^n,\ x\neq 0$
$$\lim_{\e\to 0} \e \int_0^{\|x\|} t^{-1+\e} f(t) dt = f(0)=1,$$
so the functions under the integral by $x$ in $u(\e)$ converge pointwise  
to $-\ln\|x\| \hat\phi(x).$ The result follows from 
the dominated convergence theorem. \qed
\bigbreak
Now recall that 
$$v(\e)=\int_{\R^n} \frac{\|x\|^\e-1}{\e}\e \left(\int_{\|x\|}^\infty t^{-1-\e} f(t)\ dt \right) \hat\phi(x) dx.$$
We have 
$$\e \int_{\|x\|}^\infty t^{-1-\e} f(t)\ dt = \e\int_0^{1/\|x\|} t^{-1+\e} f(1/t) dt.$$
The difficulty is that we cannot apply Lemma \ref{eps} to compute the limit of
the right-hand side of the latter equality, because the function $f(1/t)$ may be discontinuous 
at zero. However, we can avoid this difficulty as follows:
\begin{lemma} \label{l3}There exist a sequence $\e_k\to 0$ and a number $c<1$ such that
$$ \lim_{k\to \infty} v(\e_k) = c \int_{\R^n}  \ln\|x\| \hat\phi(x)\ dx.$$
\end{lemma}

\pf By a dominated convergence argument, similar to the one used in the previous lemma, 
it is enough to prove that there exist a sequence $\e_k\to 0$ and a number $c<1$ such that
for every $x\in \R^n,\ x\neq 0$
$$ \lim_{k\to \infty} \e_k \int_{\|x\|}^\infty t^{-1-\e_k}f(t)dt = c.$$
For every $x\neq 0$ we have
$$\Big|\e \int_{1/\e}^{\|x\|} t^{-1-\e} f(t)\ dt\Big| \le \big|\|x\|^{-\e}-\e^\e \big| \to 0,\quad {\rm as} \quad \e \to 0,$$
so it is enough to find a sequence $\e_k$ and a number $c<1$ such that 
$$\lim_{k\to \infty} \psi(\e_k) =  c<1,$$
where
$$\psi(\e) =  \e \int_{1/\e}^\infty t^{-1-\e}f(t)dt.$$
Since the function $\psi$ is bounded by 1 on $(0,1),$
it suffices to prove that $\psi(\e)$  cannot converge to 1, as $\e\to 0.$

Suppose that, to the contrary, $\lim_{\e\to 0} \psi(\e) = 1.$ We use the following result from [VTC, p. 205]:
if $\mu$ is a probability measure on $\R^n$ and $\gamma$ is the standard Gaussian measure
on $\R^n,$ then for every $t>0$
\begin{equation} \label{vtc}
\mu\{x\in \R^n:\ |x|_2>1/t\} \le 3\int_{\R^n} \left(1-\hat\mu(ty)\right) d\gamma(y),
\end{equation}
where $|\cdot|_2$ is the Euclidean norm on $\R^n.$
Let $\mu$ be the measure satisfying $\hat\mu = f(\|\cdot\|).$ For every $\e\in (0,1),$ integrating (\ref{vtc}) we get
$$\e\int_{1/\e}^\infty t^{-1-\e} \mu\{x\in \R^n:\ |x|_2>1/t\} dt$$
\begin{equation} \label{ineq}
 \le  \int_{\R^n} \left(\e \int_{1/\e}^\infty t^{-1-\e}(1- f(t\|y\|))dt \right) d\gamma(y).
 \end{equation}
Now
$$\e\int_{1/\e}^\infty t^{-1-\e} \mu\{x\in \R^n:\ |x|_2>1/t\}\ dt$$
$$= \e\int_{0}^\e  t^{-1+\e} \mu\{x\in \R^n:\ |x|_2>t\}\ dt.$$
and, by Lemma \ref{eps}, the limit of the left-hand side of (\ref{ineq}) as $\e\to 0$ is equal to 
$\mu(\R^n\setminus\{0\}).$ 

On the other hand, the functions
\begin{equation}\label{func}
h_\e(y) = \e \int_{1/\e}^\infty t^{-1-\e}(1- f(t\|y\|))dt 
\end{equation}
are uniformly (with respect to $\e$) bounded by 2.  Write these functions as
$$ h_\e(y)= \e \int_{1/\e}^\infty t^{-1-\e}(1- f(t\|y\|))dt = \e^\e -\|y\|^{\e} \e \int_{\|y\|/\e}^\infty t^{-1-\e}f(t) dt $$
$$=  \e^\e - (\|y\|^\e-1)\e \int_{\|y\|/\e}^\infty t^{-1-\e}f(t) dt$$
$$ - \e \int_{\|y\|/\e}^{1/\e} t^{-1-\e} f(t) dt -\e\int_{1/\e}^\infty t^{-1-\e} f(t) dt.$$
For every $y\neq 0$
$$\Big|(\|y\|^{\e}-1) \e \int_{\|y\|/\e}^\infty t^{-1-\e}f(t) dt \Big|\le \big|\|y\|^\e-1\big| \left(\frac{\|y\|}{\e}\right)^{-\e}\to 0,\quad {\rm as}\quad \e\to 0,$$
$$\Big|\e \int_{\|y\|/\e}^{1/\e} t^{-1-\e} f(t) dt\Big|\le \big|\e^\e-(\|y\|/\e)^{-\e}\big|\to 0,\quad {\rm as}\quad \e\to 0,$$
and by our assumption
$$\e\int_{1/\e}^\infty t^{-1-\e} f(t) dt = \psi(\e)\to 1, \quad {\rm as}\quad \e\to 0.$$
Therefore, the functions  $h_\e$ converge to zero pointwise as $\e\to 0$ and are uniformly bounded by a constant.
By the  dominated convergence theorem, the limit of the right-hand side of (\ref{ineq}) is equal to 0, as $\e\to 0.$  

Sending $\e\to 0$ in (\ref{ineq}), we get $\mu(\R^n\setminus\{0\}) = 0,$
therefore the probability measure $\mu$ is a unit atom at the origin and $f$ is a constant function,
which contradicts to the assumption of Theorem \ref{main}. \qed
\bigbreak
{\bf End of the proof of Theorem \ref{main}:}  Let $\e_k$ be the sequence from
Lemma \ref{l3}. Recall that $g$ is a non-negative function (see (\ref{nonneg})).  By Lemmas \ref{l1}, \ref{l2}, \ref{l3},
$$0\le \lim_{k\to \infty} g(\e_k) =\lim_{k\to \infty} (u+v+w)(\e_k) = (-1+c)\int_{\R^n} \ln\|x\| \hat\phi(x) dx,$$
where $c<1.$ Therefore,
$$ \big\langle (\ln\|\cdot\|))^\wedge, \phi \big\rangle = \int_{\R^n} \ln\|x\| \hat\phi(x) dx  \le 0$$
for every even non-negative test function $\phi$ supported outside of the origin.
By Proposition \ref{log}, $(\R^n,\|\cdot\|)$ embeds in $L_0.$ \qed
\bigbreak

\section{Examples}

The concept of embedding of a normed space in $L_0$ was studied in [KKYY].  In particular,
it was proved in [KKYY, Th.6.7] that

\begin{pr} \label{q>0} Every finite dimensional subspace of $L_p,\ 0<p\le 2$ embeds in $L_0.$ 
\end{pr}

On the other hand, as proved in [KKYY,Th.6.3],

\begin{pr} \label{q<0} If $(\R^n,\|\cdot\|_K)$ embeds in $L_0$, it also embeds in $L_p$ for every
$-n<p<0.$ 
\end{pr}
The definition and properties of embeddings in $L_p,\ p<0$ and their connections with
geometry can be found [K3, Ch. 6]. 
Propositions \ref{q>0} and \ref{q<0}  confirm the place 
of $L_0$ in the scale of $L_p$-spaces. Speaking informally, the space $L_0$ is larger than 
every $L_p,\ p\in (0,2)$, but smaller than every $L_p,\ p<0.$ 

There are many examples of normed spaces that embed in $L_0,$ 
but don't embed in $L_p,\ p\in (0,2)$ (see [KKYY, Th. 6.5]). In particular, the spaces
$\ell_q^3,\  q>2$ have this property. In fact, every three dimensional normed space
embeds in $L_0$ (see [KKYY, Corollary 4.3]). However, starting from
dimension 4, there are many normed spaces that do not embed in $L_0.$
The following result from [K3, Th. 4.19] essentially shows that a normed space 
with dimension greater than 4 does not embed in $L_0$ if the second derivative
of its norm at zero in at least one direction is equal to 0. 
 
\begin{pr} \label{sdt}
Let $n\ge 4,\ -n<p<0$ and let $X=(\R^n,\|\cdot\|)$ be an $n$-dimensional  normed
space with a normalized basis $e_1, \dots, e_n$ so that: 
\item{(i)} For every fixed
$(x_2,\dots,x_n)\in \R^{n-1}\setminus\{0\},$ the 
function 
$$x_1\mapsto \|x_1e_1+\sum_{i=2}^n x_ie_i\|$$ 
has a continuous second derivative everywhere on $\R,$ and 
$$\|x\|_{x_1}^{'}(0,x_2,\dots,x_n) =
\|x\|_{x_1^2}^{''}(0,x_2,\dots,x_n)= 0,$$ 
where
$\|x\|_{x_1}^{'}$ and $\|x\|_{x_1^2}^{''}$  stand for the
first and second partial derivatives by $x_1$ of the norm
$\|x_1e_1+\dots+x_ne_n\|.$ 

\item{(ii)} There exists a constant $C$ so that,
for every $x_1\in \R$ and every $(x_2,\dots,x_n)\in \R^{n-1}$ with
$\|x_2e_2+\dots+x_ne_n\|=1,$ one has
$$\|x\|_{x_1^2}^{''}(x_1,x_2,\dots,x_n) \le C .$$

\item{(iii)} Convergence in the limit 
$$\lim_{x_1\to 0} \|x\|_{x_1^2}^{''}
(x_1,x_2,\dots,x_n)= 0$$ 
is uniform with respect to $(x_2,\dots,x_n)\in \R^{n-1}$
satisfying the condition $\|x_2e_2+\dots+x_ne_n\|=1.$ 
\medbreak
Then the space $(\R^n,\|\cdot\|)$ does not embed in $L_0.$ 
\end{pr}

\pf It was proved in [K3, Th. 4.19] that under the assumptions of Proposition \ref{sdt}
the function $\|\cdot\|_K^{-p}$ represents a positive definite distribution if and
only if $p\in (n-3,n].$ By [K3, Th. 6.15] the space $(\R^n,\|\cdot\|_K)$ does not embed 
in $L_p,\ p\in (-1,0),$ so it also does not embed in $L_0$ by Proposition \ref{q<0}. 
The result follows from Theorem \ref{main}. \qed
\bigbreak

From Proposition \ref{sdt} and Theorem \ref{main} we immediately get
\begin{co} \label{n>3} If a normed space $(\R^n,\|\cdot\|),\ n\ge 4$ satisfies the conditions
of Proposition \ref{sdt}, then a function of the form $f(\|\cdot\|)$ can be positive definite 
only if $f$ is a constant function. The norm of such a space cannot appear as
the standard of an $n$-dimensional version.
\end{co}

\bigbreak

Let us give several examples of spaces satisfying the conditions of 
Proposition \ref{sdt}.  For normed spaces $X$ and $Y$ and 
$q\in \R,\ q\ge 1,$ the $q$-sum $(X\oplus Y)_q$ of $X$ 
and $Y$
is defined as the space of pairs $\{(x,y):\ x\in X, y\in Y\}$ with the norm
$$\|(x,y)\|= \left(\|x\|_X^q+\|y\|_Y^q \right)^{1/q}.$$  It was proved in [K2, Th 2]
that such spaces with $q>2$ satisfy the conditions of Proposition \ref{sdt}
provided that the dimension of $X$ is greater or equal to 3.

Another example is that of Orlicz spaces. Recall that an {\it Orlicz function} $M$ is a non-decreasing 
convex function on $[0,\infty)$ such that $M(0)=0$ and 
$M(t)>0$ for every $t>0.$
The norm $\|\cdot\|_M$ of the $n$-dimensional
Orlicz space $\ell_M^n$ is defined implicitly by the equality
$$\sum_{k=1}^n M(|x_k|/\|x\|_M) = 1,\ x\in \R^n\setminus \{0\}.$$
As shown in [K2, Th 3], the spaces $\ell_M^n,\ n\ge 4$ 
satisfy the conditions of Proposition \ref{sdt} if the Orlicz function
$M\in C^2([0,\infty))$ is such that $\ M'(0) = M''(0) =0.$ 

\begin{co} If a normed space $(\R^n,\|\cdot\|)$ contains a subspace isometric to 
$(X\oplus Y)_q$, where $q>2$ and the dimension of $X$ is at least 3, or contains
an Orlicz space $\ell_M^4,$ where $M$ is an Orlicz function such that $M\in C^2([0,\infty))$ 
and $\ M'(0) = M''(0) =0,$ then a function of the form $f(\|\cdot\|)$ can be positive definite
only if $f$ is a constant function. 
\end{co}


\begin{thebibliography}{99}

\bibitem[AMM]{AMM} { I.~Aharoni, B.~ Maurey and B.~Mityagin},
\textit{Uniform embeddings of metric spaces and of Banach 
spaces into Hilbert spaces}, Israel J. Math. \textbf{52} (1985), 251--265. 

\bibitem[BDK]{BDK} { J.~Bretagnolle, D.~Dacunha-Castelle and J.~L.~Krivine},
\textit{Lois stables et espaces $L_p$},
Ann. Inst. H. Poincar\'e  Probab.  Statist.  \textbf{2} (1966),  
231--259. 

\bibitem[CKS]{CKS} { S.~Cambanis, R.~Keener, and G.~Simons},
\textit{On $\alpha$-symmetric multivariate distributions},
J. Multivariate Analysis \textbf{13} (1983),  213--233.

\bibitem[D]{D} { L.~Dor},
\textit{Potentials and isometric embeddings in $L_{1}$},
Israel J. Math. \textbf{24} (1976), 260--268.

\bibitem[E]{E} { M.~Eaton},
\textit{On the projections of isotropic distributions},
Ann. Stat. \textbf{9} (1981), 391--400. 

\bibitem[F]{F} { T.~S.~Ferguson}, \textit{A representation of
the symmetric bivariate Cauchy distributions},
Ann. Math. Stat. \textbf{33} (1962), 1256--1266.

\bibitem[G]{G} { T.~Gneiting}, \textit{On $\alpha$-symmetric multivariate 
characteristic functions},  J. Multivariate Anal. \textbf{64}  (1998), 131--147.

\bibitem[H]{H} { C.~Herz,} \textit{A class of negative definite functions},
Proc. Amer. Math. Soc. {\bf 14}  (1963), 670--676.

\bibitem[KKYY]{KKYY} { N.~J.~Kalton, A.~Koldobsky, V.~Yaskin and M.~Yaskina,}
\textit{The geometry of $L_0$}, Canad. J. Math. {\bf 59} (2007), 1029--1049.

\bibitem[K1]{K1} { A.~ Koldobsky},  \textit{ The  Schoenberg problem on 
positive-definite functions}, Algebra i Analiz \textbf{ 3} (1991), 78--85; 
translation in St. Petersburg Math. J. \textbf{3} (1992), 563--570. 

\bibitem[K2]{K2} { A. Koldobsky}, \textit{Second derivative test for intersection 
bodies}, Advances in Math. \textbf{136} (1998), 15--25.

\bibitem[K3]{K3} {A.~Koldobsky}, \textit{Fourier analysis in convex geometry},
Amer. Math. Soc., Providence RI, 2005.

\bibitem[K4]{K4} { A. Koldobsky}, \textit{A note on positive definite norm dependent functions}, 
Proceedings of the Conference on High Dimensional Probability, Luminy, 2008, to appear

\bibitem[KL]{KL} { A.~Koldobsky and Y.~Lonke},  \textit{A short proof of
Schoenberg's conjecture on positive definite functions},
Bull. London Math. Soc.  \textbf{31}  (1999), 693--699.

\bibitem[Ku]{Ku} { Yu.~G.~Kuritsyn},
\textit{Multidimensional versions and two problems of Schoenberg},
Problems of Stability of Stochastic Models, VNIISI, Moscow, 1989,
72--79.  

\bibitem[Le]{Le} { P.~L\'evy}, \textit{Th$\acute {e}$orie de l'addition de variable 
al$\acute {e}$atoires}, Gauthier-Villars, Paris, 1937.

\bibitem[L]{Li} { J.~Lindenstrauss,}
{\em On the extension of operators with finite dimensional range},
Illinois J. Math. {\bf 8}  (1964), 488--499.

\bibitem[Li1]{Li1} { A.~Lisitsky}, 
\textit{One more proof of Schoenberg's conjecture},
unpublished manuscript, 1991.

\bibitem[Li2]{Li2} { A.~Lisitsky}, 
\textit{The Eaton problem and multiplicative properties of multivariate distributions},
Theor. Probab. Appl.  \textbf{42} (1997), 618--632.

\bibitem[M1]{M1}  { J.~Misiewicz}, 
\textit{On norm dependent positive definite functions},
Bull. Acad. Sci. Georgian SSR \textbf{130} (1988), 253--256.

\bibitem[M2]{M2} { J.~Misiewicz}, 
\textit{Positive definite functions on $\ell_{\infty}$},
Stat. Probab. Let. \textbf{8} (1989), 255--260.

\bibitem[M3]{M3} { J.~Misiewicz}, \textit{Substable and pseudo-isotropic 
processes---connections with the geometry of subspaces of $L\sb \alpha$-spaces}, 
Dissertationes Math. (Rozprawy Mat.)  \textbf{358}  (1996).

\bibitem[R]{R} { D.~St.~P.~Richards}, 
\textit{Positive definite symmetric functions on finite 
dimensional spaces. 1.Applications of the Radon transform},
J. Multivariate Analysis \textbf{19} (1986),  280--298. 

\bibitem[S1]{S1} { I.~J.~Schoenberg},
\textit{Metric spaces and positive definite functions},
Trans. Amer. Math. Soc. \textbf{44} (1938), 522--536.

\bibitem[S2]{S2}  { I.~J.~Schoenberg},
\textit{Metric spaces and completely monotone functions},
Annals of Math. \textbf{39} (1938), 811--841.

\bibitem[VTC]{VTC} { N.~N.~Vakhania, V.~I.~Tarieladze and S.~A.~Chobanyan}, \textit{Probability
distributions on Banach spaces}, D.~Reidel Publishing Company, Dordrecht, 1987.

\bibitem[Z1]{Z1} { V.~Zastavnyi}, 
\textit{Positive definite norm dependent functions},
Dokl. Russian Acad. Nauk.  \textbf{325} (1992), 901--903. 

\bibitem[Z2]{Z2} { V.~Zastavnyi}, \textit{Positive definite functions depending 
on the norm}, Russian J. Math. Phys. \textbf{1} (1993), 511--522.


\end{thebibliography}
\end{document}